\documentclass[letterpaper, 10 pt, conference]{ieeeconf}  

\IEEEoverridecommandlockouts                              

\usepackage{graphicx} 
\usepackage{amsmath} 
\usepackage{amssymb}  
\usepackage{color}
\usepackage{standalone}
\usepackage{pgfplots}
\usepackage{tikz}
\usepackage{booktabs}
\usepackage{multirow}
\usepackage{todonotes}
\usetikzlibrary{shapes,arrows}
\usetikzlibrary{fit,positioning}
\pgfplotsset{compat=newest}

\title{\LARGE \bf
Optimization Strategies for Real-Time Control of an Autonomous Melting Probe 
}

\author{Christian Meerpohl, Kathrin Fla\ss kamp and Christof B\"uskens
\thanks{C. Meerpohl, K. Fla\ss kamp and C. B\"uskens are with the Center for Industrial Mathematics, 
        University of Bremen, 28359 Bremen, Germany,
        {\tt\small \{cmeerpohl,bueskens\}@math.uni-bremen.de},
        {\tt\small kathrin.flasskamp@uni-bremen.de}}%
}

\newcommand\copyrighttext{%
	\footnotesize \textcopyright \,2018 IEEE. Personal use of this material is permitted. Permission from IEEE must be obtained for all other uses, in any current or future media, including reprinting/republishing this material for advertising or promotional purposes, creating new collective works, for resale or redistribution to servers or lists, or reuse of any copyrighted component of this work in other works.}
\newcommand\copyrightnotice{%
	\begin{tikzpicture}[remember picture,overlay]
	\node[anchor=north,yshift=-10pt] at (current page.north) {\fbox{\parbox{\dimexpr\textwidth-\fboxsep-\fboxrule\relax}{\copyrighttext}}};
	\end{tikzpicture}%
}

\begin{document}

\maketitle
\copyrightnotice
\vspace*{-10pt}

\begin{abstract}

We present an optimization-based approach for trajectory planning and control of a maneuverable melting probe with a high number of binary control variables. The dynamics of the system are modeled by a set of ordinary differential equations with a priori knowledge of system parameters of the melting process. The original planning problem is handled as an optimal control problem. Then, optimal control is used for reference trajectory planning as well as in an MPC-like algorithm. Finally, to determine binary control variables, a MINLP fitting approach is presented. The proposed strategy has recently been tested during experiments on the Langenferner glacier. The data obtained is used for model improvement by means of automated parameter identification.
\end{abstract}

\section{Introduction}
\label{sec:introduction}
During the Cassini-Huygens mission, which lasted for nearly 20 years and came to an end by September 2017, the Saturnian moon Enceladus has been found to be one of the most interesting celestial bodies in our solar system regarding future exploration missions.
Plumes that ascend from the south polar region could be detected by the Cassini spacecraft and were identified to contain ionized water vapor. Additionally, camera recordings showed large cracks in Enceladus' kilometer thick icy shell with higher surface temperatures, the so-called ``Tiger stripes''. 
An evaluation of the gathered information indicates cryovolcanism as the reason for such activity, leading to the conclusion that there likely is a global liquid water ocean below the surface \cite{porco2006cassini}. The possibility of proving extraterrestrial life or habitability was a starting signal for several subsequent mission concepts.\par
\begin{figure}[ht]
	\centering
	\includegraphics[width=5.5cm]{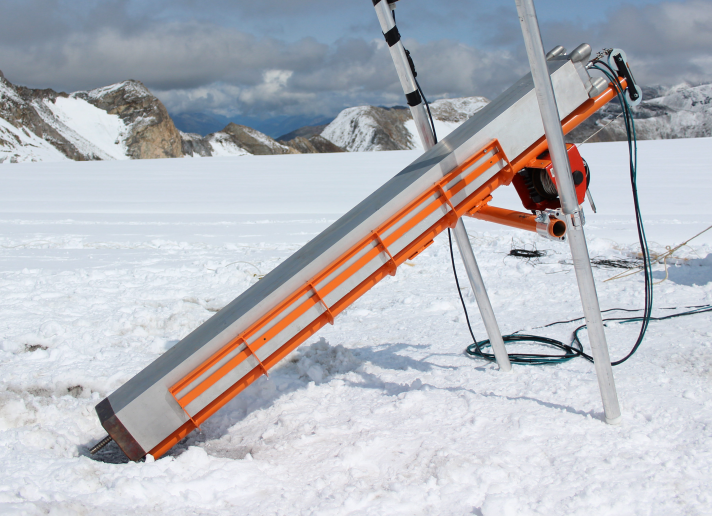}
	\caption[icemole]{Melting probe IceMole in a launchpad on the Langenferner Glacier in Italy during the 2017 test run. Copyrights: Joachim Clemens, University of Bremen.}
	\label{fig:icemole}
\end{figure}

 
One of these ideas has been developed during the Enceladus Explorer (``EnEx'') collaborative project. A possible scenario would be to retrieve water samples from the ocean by drilling through the icy shell. If aiming for cracks close to the geysers and hitting liquid water reservoirs, the mission objective could be accomplished relatively easily. 
Therefore, a maneuverable melting probe has been built, the IceMole \cite{dachwald2014icemole} (Fig. \ref{fig:icemole}).
During several field tests, its ability to locate itself and to maneuver through ice has been shown. A scientific highlight was the first successful retrieving of a clean brine sample of the Blood Falls in Antarctica \cite{kowalski2016navigation}. 
After demonstrating the technical feasibility of the mission concept in principle, necessary technologies shall now be further enhanced in six follow-up projects, each focusing on different parts.
In EnEx-CAUSE, a fully autonomous navigation system is developed, involving sensor fusion, decision making and steering of the probe. 
\par 
Guiding the IceMole raises a sequence of challenging tasks. From modeling the system dynamics, to trajectory planning, up to determining the actuator controls, a coherent concept is required. A proven strategy to handle the problems of trajectory planning and model predictive control is given by the definition and solution of optimal control problems (OCP). Nowadays, complex and highly nonlinear systems do not necessarily represent an obstacle against the finding of feasible and even optimal solutions under real-time conditions. The combination of a binary character and the amount of 32 control variables, however, makes the IceMole a challenging candidate for real-time control. As a consequence, a methodology to deal with potentially high-dimensional mixed-integer control problems (MIOCP), that have to be solved within time restrictions, has to be devised.
\subsection{Related Work}
The study of OCPs became popular in the 1950s and 1960s, when there was a high demand for solutions in the military sector \cite{sussmann1997300}. Solution strategies can be divided into two classes of approaches \cite{stryk1992direct}: On the one hand, direct methods transcribe the problem by collocation, multiple shooting or full discretization methods into nonlinear programs (NLP). These potentially high-dimensional problems can be solved by sequential quadratic programming (SQP) or interior point algorithms. On the other hand, the indirect approach uses Pontryagin's maximum principle to state the necessary optimality conditions, which lead to a boundary value problem that again can be solved by shooting methods. MIOCPs, often also referred to as hybrid optimal control problems, have firstly been studied in the 1980s. They are considered to be much harder to solve than pure OCPs \cite{kirches2011fast}.
For solving them, there is no universal best approach but a variety of strategies. By analogy with OCPs, they can be discretized to mixed-integer nonlinear programs (MINLP), which typically are solved by branch-and-bound algorithms. Alternatively, there are strategies such as relaxation, transformation or rounding techniques and combinations thereof, as e.g. proposed in \cite{gerdts2006variable} and \cite{sager2005numerical}.
\subsection{Contribution}
Our contribution is twofold. Firstly, we present a methodology on how the IceMole melting probe can be autonomously steered during glacier tests using a combination of off-the-shelf optimization software. 
We propose a pragmatic concept that at some point neglects dynamics exactness, but guarantees feasible solutions under real-time conditions. By relaxing a MIOCP, we are able to plan and adjust optimal trajectories in a robust and efficient way. The binary variables are determined in a subsequent fitting algorithm. We set up a MINLP that is being frequently solved so that feasible and effective binary control configurations can be applied.\par Secondly, we conduct a post-test analysis. Because of the limited repeatability of the test procedure, we do a comprehensive evaluation of the obtained data in order to enhance our algorithms for future field tests. On the one hand, this includes revealing and adjustment of weak aspects coming with our methodology. On the other hand, we perform automated parameter identification to improve our system dynamics model.
\section{System Dynamics} 
\label{sec:dynamics}
\subsection{Melting Probes}
Melting probes are designed to reach a certain destination in an icy environment in order to perform a particular action. Typically, they allow only straight-downward motion.  This type of movement has been studied in \cite{aamot1967heat} and \cite{ulamec2007access} with respect to the correlation between induced heating power and melting velocity under the effects of conductive loss. In order to actively maneuver through ice, a concept around curvilinear motion is required. For instance, there are studies about the use of heating plates. One \cite{kumano2005direct} deals with the application of an asymmetrically distributed force, another  \cite{schuller2016curvilinear} presents a framework around differential heating and a constant contact force. By applying a temperature gradient at the surface of a pressurized plate, a curvilinear melting path can be enforced, which can be described by the curve radius and the mean velocity. However, transferring these principles to a universal model for trajectory planning of a specific melting probe in three dimensions is not obvious. The inertia of the body cannot be disregarded, as well as specific attributes that come along with the IceMole's construction. 
\subsection{IceMole Model and Constraints}
The IceMole is a special melting probe, as it uses two complementary advance strategies. A screw guarantees a steady contact force towards the ice front and, in principle, allows movement against the weight, while a collection of heater elements determines the actual speed and, additionally, allows curvilinear melting due to the differential heating concept. 
\begin{figure}[ht]
	\centering
	\hspace*{.1cm}
	\resizebox{90mm}{!}{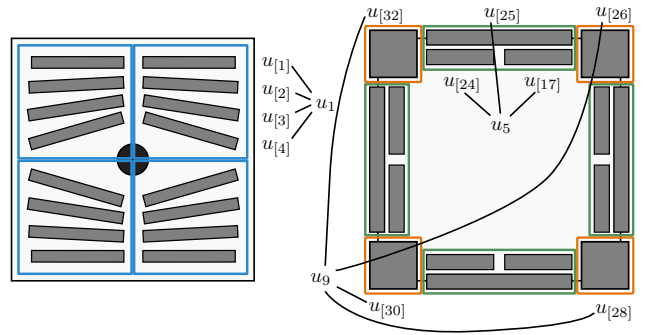}
	\caption{Schematic diagram of the IceMole's layout and heater numbering. \textit{Left}: Screw in the center and all 16 head heater units including an exemplary grouping. \textit{Right}: Wall heater elements on the inner side (two of each), while the corresponding back heater elements are illustrated on the outer side. Four units are positioned in the corners.}
	\label{fig:icemole_scheme}
\end{figure} 
\par
In Fig. \ref{fig:icemole_scheme}, the IceMole's layout, including the positioning of the screw and all heater elements $\mathcal{I}_T=\{1,\dots,32\}$, is illustrated. 
The heater units are numbered clockwise, beginning with those disposed in the head at the front side 
$\mathcal{I}_H=\{1,\dots,16\}$, followed by those positioned in the walls $\mathcal{I}_W= \{17,\dots,24\}$ and completed by those in the back part $\mathcal{I}_B=\{25,\dots,32\}$.
Before we model the dynamic behavior in an efficient way, we make a simplification. The head heater units in each quadrant are grouped together, so that there are four groups $\mathcal{I}_{1},\dots,\mathcal{I}_{4}$. We attach one control variable $u_i, i=1,\dots,4$ to each of these groups. The same can be done with the wall heater units and the corresponding back heater unit on each side. Since the units in the back corners stick out of the body, they constantly have to be turned on ($u_9 \equiv 1$). Otherwise, the IceMole could get stuck or, at least, be restricted in its movement. Therefore, they are not considered as control variables for modeling purposes. Altogether, we define eight control variables $\boldsymbol{u}=\{u_i\}_{i=1,\dots,8}$. \par
To model the system dynamics, we do not want to examine the melting process itself, involving concepts of thermal conductivity, convection and the transfer of energy by phase changes. Instead, we will set up a relation between the observable state of the system and the actuator controls by using a priori knowledge of the system. 
Assuming that the IceMole's translation can only occur along its longitudinal axis and neglecting gravitation, we describe our motion model by the kinematic equations of motion with six degrees of freedom in an east-north-up frame \cite{stengel2015flight} as: 
\begin{equation}
\label{eq:ode}
\frac{d}{dt}
\begin{pmatrix}
   x\\ y\\ z\\ q_w\\ q_x\\ q_y\\ q_z
\end{pmatrix} = 
\begin{pmatrix}
&(1-2q_y^2-2q_z^2)v \\
& 2(q_xq_y+q_wq_z)v \\
& 2(q_xq_z-q_wq_y)v \\
& \frac{1}{2}(-q_x\omega_x-q_y\omega_y-q_z\omega_z) \\
& \frac{1}{2}(q_w\omega_x+q_y\omega_z-q_z\omega_y) \\
& \frac{1}{2}(q_w\omega_y+q_z\omega_x-q_x\omega_z) \\
& \frac{1}{2}(q_w\omega_z+q_x\omega_y-q_y\omega_x) \\
\end{pmatrix}
\end{equation}
It is stated as an autonomous ordinary differential equation system ${\dot{\boldsymbol{x}} = f(\boldsymbol{x},\boldsymbol{u},\boldsymbol{p})}$.
The state vector $\boldsymbol{x}$ is a pose, which is defined by coordinates $x,y,z$ and a quaternion $\vec q = (q_w,q_x,q_y,q_z)$. Quaternions allow a non-unique but differentiable representation of rotations in three dimensions. \par
It is still unclear how the heater controls are related to the kinematic equations.
There is a general understanding of how each unit affects the system's behavior. Namely, the IceMole's speed is mainly defined by the front heaters, while its curvilinear motion is additionally determined by the wall heaters. The heaters on the back panel support backwards motion when the probe has to be pulled out of a partially frozen melting channel. During the EnEx collaborative project, a proven control strategy was to either operate on all head heaters for straight forward movement or to turn on half of the head heaters and the corresponding wall heaters for curvilinear melting. 
Since a nearly linear dependence between the melting velocity and the induced heat has been pointed out in \cite{aamot1967heat}, we model this relationship by
\begin{equation}
\label{eq:velocity}
v = \sum\limits_{i \in \{1,\dots,4\}} u_i\frac{v_{\max}}{4},\quad 0\le u_i\le 1,
\end{equation}
where $v_{\max}$ is the maximal velocity. \\
Due to the torque applied by the screw, the IceMole holds a certain roll rate. Its extent depends on the ice quality, on the one hand, and on whether wall heaters are used or not, on the other hand. As previous tests have shown, those effects can change locally and cannot be characterized from data validation at this point. Therefore, we make a very simple generalization of a linear increase with respect to the velocity:
\begin{equation}
\label{eq:rollrate}
\omega_x = v \cdot \eta, \quad \eta > 0.
\end{equation}
Assuming that we know the radius for maximal curvilinear melting $R_{\min}$ and that rotation can only occur during forward movement, we define the angular velocities with respect to the y- and z-axis by
\begin{equation}
\label{eq:pitchrate}
\omega_y = \frac{v}{R_{\min}} \frac{\alpha}{2}({u_1+u_4-u_2-u_3})+(1-\alpha)(u_5-u_7),
\end{equation}
\begin{equation}
\label{eq:yawrate}
\omega_z = \frac{v}{R_{\min}} \frac{\alpha}{2}({u_1+u_2-u_3-u_4})+(1-\alpha)(u_6-u_8).
\end{equation}
The parameter $\alpha \in [0,1]$ describes to which extent the curvilinear melting process can be attributed to the head and wall heater units respectively. In total, our model includes four parameters ${\boldsymbol{p}=\{v_{\max},\eta,R_{\min},\alpha\}}$. \par
Power supply is restricted so that not all units can be operated at once. 
Therefore, power consumption for each group has to be incorporated, meaning that the weighted sum of the control variables has to be restricted:
\begin{equation}
\label{eq:powerallgroup}
\sum\limits_{i = \{1,\dots,8\}}P_{i} u_{i}(t) \le \tilde P_{T,\max},
\end{equation}
where $P_i = \sum\nolimits_{j \in \mathcal{I}_i} P_{[j]}$ is the sum of the single units related to each group and $\tilde P_{T,\max}$ is the maximal total power consumption $P_{T,\max}$ minus the power, which has to be reserved for the back heater corner units $P_9$.
Additionally, the IceMole's software interface demands that power has to be separately reserved for the wall and back heaters. This leads to one further constraint:
\begin{equation}
\label{eq:powerwallgroup}
\sum_{i=\{5,\dots,8\}}^{} P_i u_i(t) \le \tilde P_{W,\max},
\end{equation}
with 
$ \tilde P_{W,\max} = P_{W,\max}+(P_{B,\max}-P_9)$.

\section{Trajectory Planning}
\label{sec:planning}
With the given dynamics, trajectories can be planned. Motion planning algorithms like PRM (Probabilistic roadmap) planners \cite{kavraki1996probabilistic} or RRTs (Rapidly-exploring random trees) \cite{lavalle98rapidly} are considered to be state of the art to formulate and solve these problems in many applications. Since we do not take any possible obstacles into account, the main difficulty consists in complying with the dynamics. The system of differential equations \eqref{eq:ode} has to be solved, as well as the constraint equations \eqref{eq:velocity} - \eqref{eq:powerwallgroup} for every time step. Additionally, one might want to plan trajectories with respect to different criteria. This favors a formulation as an optimization problem, allowing to quickly solve all at once. \par
Until this point, we have not looked closely at the possible control variables. The heater elements that are plugged into the IceMole can either be switched on or off, which means they should be modeled as binary controls. Under this assumption, the planning problem has to be stated as a MIOCP. Compared to OCPs, this lifts the planning problem to another level. Continuous optimal control problems can be solved for millions of variables and constraints in less than a second, while the computation time for MINLPs in principal grows exponentially with every additional optimization variable \cite{bixby2004mixed}. Thus, we relax the MIOCP. The continuous heater controls will further be referred to as ${\boldsymbol{u} \in [0,1]^8}$, while the binary controls are stated as ${\boldsymbol{v} \in \{0,1\}^8}$. 
An OCP is defined in its standard form by
\begin{equation}\tag{OCP}\label{eq:OCP}
\begin{array}{lrl}
\min\limits_{x,u,t_f}  & F(x(t),u(t),t)\\
s.t. & 	 \dot{x}(t) & = f(x(t),u(t),t),\\
&	\omega(x(0),x(t_f)) &= 0, \\
& \psi(x(t),u(t))&\le 0,\\ \displaystyle
&  \multicolumn{2}{l}{x \in \mathbb{R}^{n_x}, u \in \mathbb{R}^{n_u}, t \in [0,t_f].}
\end{array}
\end{equation} 
To solve a problem of this form, we use discretization methods, e.g. the trapezoidal rule to transcribe it into a NLP of general form
\begin{equation}\tag{NLP}\label{eq:NLP}
\begin{array}{llll}
	\min\limits_{z \in \mathbb{R}^{n_z}} & f(z) \\ 
	s.t. & g_i(z) &\le 0, & i=1 ,..., m_g,\\
	     & h_j(z) &= 0, &j=1 ,..., m_h.
\end{array}
\end{equation}
Common techniques to solve these sparse and potentially high-dimensional NLPs are sequential quadratic programming or interior point algorithms. We use the optimization routine WORHP \cite{buskens2012esa} with an SQP method that uses an interior point algorithm on the QP level. 
Below, we will deduce our trajectory planning problem. We distinguish between two formulations in which one version is more restrictive and supposed to be solved only once, while the other one is executed frequently.
\subsection{Reference Trajectory}
Reaching a fixed goal pose $\boldsymbol{x}_f$ with an exact orientation is nearly impossible because of the IceMole's limited controllability. Additionally, the final roll angle $\phi$ is of no particular interest. Therefore, only pitch and yaw angle $\theta, \psi$ shall be considered. The attitude deviation in the final state can be interpreted as a rotation. We use the boxminus operator ${\boxminus:SO(3)\times SO(3)\to \mathbb{R}^3}$ (\cite{bloesch2016primer}) to determine its representation in Euler angles: 
\begin{equation}
\label{eq:boxminus}
(\Delta\phi_f,\Delta\theta_f,\Delta\psi_f) := \vec q(t_f) \boxminus \vec q_f,
\end{equation}
where $\Delta\phi_f,\Delta\theta_f,\Delta\psi_f$ are the deviations with respect to the three Euler angles. 
Now, $\Delta\theta_f,\Delta\psi_f$ can be restricted to a certain value $\varphi_{max}$:
\begin{equation}
\label{eq:finalattitude}
-\varphi_{max} \le \Delta\theta_f,\Delta\psi_f \le \varphi_{max}
\end{equation}
Additionally, we introduce box constraints for the translational velocity and the angular velocities:
\begin{equation}
\label{eq:velocityconstraints}
\begin{array}{rcl}
0 \le& v &\le v_{max} \\
-\omega_{i,max} \le& \omega_{i} &\le \omega_{i,max}\quad \forall i \in \{x,y,z\} 
\end{array}
\end{equation}
We define the values $v,\omega_{x},\omega_{y},\omega_{z}$ as control variables, since an increase in the number of linear equations is often numerically preferable over fewer nonlinear equations. Still unknown is the objective of the OCP. 
Since a linear dependence between the mean melting velocity and power consumption has been modeled, a compelling energy-optimal solution can hardly be obtained with respect to this model. Therefore, we choose a time-optimal formulation and add a weighted penalty term (with ${\beta > 0}$) for the final orientation deviation. Altogether, with a state vector ${\boldsymbol{x}(t) := (x,y,z,q_w,q_x,q_y,q_z)(t)}$ and a control vector ${\tilde{\boldsymbol{u}}(t) := (v,\omega_x,\omega_y,\omega_z,\boldsymbol{u})(t)}$, we define our trajectory planning problem by 
\begin{equation}\label{eq:tpp}\tag{TPP}
\begin{array}{lll}
\min\limits_{\boldsymbol{x},\tilde{\boldsymbol{u}},t_f}  &
t_f + \beta({\Delta\theta_f}^2 + {\Delta\psi_f}^2) \\
s.t. & \text{eq. } \eqref{eq:ode} -\eqref{eq:velocityconstraints} \text{ hold}\, \forall t \in [0,t_f], \\ &
\boldsymbol{x}(0) = \boldsymbol{x}_0, \\
& (x,y,z)(t_f) = (x_f,y_f,z_f).
\end{array}
\end{equation}
\subsection{Replanning}
To follow a computed path, a replanning or feedback control algorithm has to be implemented. If obstacles are considered, trajectory updates can involve large discrepancies. Since we neglect the occurrence of obstacles, the potential for error only arises from the dynamics, control application and environmental conditions. Therefore, we can expect a new solution to lie close enough to the previous one, which leads to a high convergence rate. However, the problem can become locally infeasible. Since the melting probe's possibilities to change direction are rather limited, reaching a fixed destination might become more and more difficult with decreasing distance. Therefore, we remove the final position constraint and, instead,
add another penalty term (with ${\gamma > 0}$) to the cost function. One could still bound this term, but we leave this decision to a higher-level decision making module.
Given a current pose update $\boldsymbol{x}_k$, the trajectory replanning problem at time step $t_k$ can be stated as
\begin{equation}
\label{eq:trpp}
\tag{MPC}
\begin{array}{lll}
\min\limits_{\boldsymbol{x},\boldsymbol{\tilde u},t_f^{[k]}} & t_f^{[k]} +\beta({\Delta\theta_f}^2+{\Delta\psi_f}^2) \\[-10pt]
& + \gamma({\Delta x_f}^2+{\Delta y_f}^2+{\Delta z_f}^2)  \\[5pt]
s.t. & \text{eq. } \eqref{eq:ode} -\eqref{eq:velocityconstraints} \text{ hold}\, \forall t \in [t_k,t_f^{[k]}], \\ &
\boldsymbol{x}(t_k) = \boldsymbol{x}_k. 
\end{array}
\end{equation}
It can be interpreted as a model predictive control algorithm with an infinite prediction and control horizon and weights not equal to zero only for the final state (cf. \cite{grune2011nonlinear}). 

\section{Actuator Control}
\label{sec:controlling}
Given the trajectories generated by \ref{eq:trpp}, a way to transform the real-valued control functions into binary control variables has to be found. This problem is closely related to pulse width modulation (PWM) techniques but with the restriction that signals cannot be separated and are connected through additional constraints. Sager \cite{sager2005numerical} proposes sum up rounding strategies for either singular, independent control functions or those who are related by SOS1 constraints ($\sum\nolimits_{i=1}^N x_i \le 1,\,x_i \in \{0,1\}$).
In our case, however, there are three inequalities of a general linear form, meaning that those strategies will not be sufficient. Therefore, we propose a MINLP optimization approach, which generates feasible control variables during every trajectory calculation cycle.\par
To determine the full set of binary controls $\{v_{[i]}\}_{i \in \mathcal{I}_T}$, we ungroup the previously defined control variables:
\begin{equation}
\label{eq:ungroup}
u_{[j]} = u_i \quad \forall j \in \mathcal{I}_i,\quad i=\{1,\dots,9\}
\end{equation}
Consequently, the power consumption constraints \eqref{eq:powerallgroup} and \eqref{eq:powerwallgroup} will be restated and complemented as follows:
\begin{equation}
\label{eq:powerall}
\begin{array}{rcl}
\sum\limits_{i \in \mathcal{I}_T}^{} P_{[i]} v_{[i]} &\le & P_{T,\max}, \\
\sum\limits_{i \in \mathcal{I}_W}^{} P_{[i]} v_{[i]} &\le & P_{W,\max}, \\
\sum\limits_{i \in \mathcal{I}_B}^{} P_{[i]} v_{[i]} &\le & P_{B,\max}.
\end{array}
\end{equation}
Although we do not consider the dynamic equations at this point, we still have to account for their physical reasoning. If we look at a short time interval of adequate length $\delta$, we demand that the defect for every control function $u(t)$, respectively $v(t)$, vanishes:
\begin{equation}
\label{eq:defect}
\int\limits_{t-\delta}^{t} (u(s))ds - \int\limits_{t-\delta}^{t} (v(s))ds \approx 0 
\end{equation}
From now on, we inspect a replanning step ${t_k \to t_{k+1}}$ of fixed length $T$.
The control function $u(t)$ is evaluated by applying a linear interpolation to the discretized solution of \ref{eq:trpp} and is extended constantly over the integral bounds.\\
If we postulate term \eqref{eq:defect} to be minimal over the whole length for all control variables and introduce weights ${\mu_{[i]} > 0}$, this leads to the following cost function:
\begin{equation} 
\label{eq:approximation}
\sum\limits_{i \in \mathcal{I}_T}^{}\mu_{[i]}^2\int_{0}^{T}\big(\int_{t-\delta}^{t} (u_{[i]}(s)-v_{[i]}(s))ds\big)^2\,dt
\end{equation}
By applying an integration scheme, e.g. the chained trapezoidal rule, we obtain a quadratic objective. 
The binary control switching points are given by the supporting points for integration. Those share the same grid for the inner and the outer integral. If constant terms are neglected, the optimization problem can be stated as 
\begin{equation}
\label{eq:bct}
\tag{BCT}
\hspace*{-.12cm}
\begin{array}{lll}
\min\limits_{v} &\sum\limits_{i \in \mathcal{I}_T}^{}\mu_{[i]}^2\sum\limits_{j=1}^{M} a^{[j]} (\sum\limits_{m=j-l}^{j} b^{[m]} (u_{[i]}^{[m]}-v_{[i]}^{[m]}) )^2  \\[10pt]
s.t. & \text{eq. }\eqref{eq:powerall} \text{ hold }\, \forall j = 1,\dots,M.
\end{array}
\end{equation}
We call it the binary control transformation (BCT).
It is solved by using the MINLP software SCIP \cite{scip2017} with WORHP on NLP level. If the computation time is too long, which is not surprising in connection with MINLPs, the algorithm can be terminated and it still provides a feasible solution.
If the replanning step size is small enough, the continuous control variables will merely change from one point to the other. Therefore, we can take the feasible solution $\{v_{[i]}^{[j]}\}^{[k]}$ and use it as an initial guess for the next calculation cycle:
\begin{equation}
\{v_{[i]}^{[j]}\}_{\text{init}}^{[k+1]} = \{v_{[i]}^{[j]}\}_{}^{[k]}
\end{equation}
As a result, we can expect solutions that become much better after a few iterations.
\section{Results}
\label{sec:results}
In August 2017, the EnEx initiative performed a field test on the Langenferner glacier in Italy, during which the IceMole was controlled fully automatically for the first time. The main objective was to show the general capability of all partners' subsystems and proper communication between them. In two ``long distance straight on'' melting tests for around 25 meters each, the correct actuating of all heater units was proven. In terms of modeling dynamics and trajectory planning, however, those tests were less informative. The only way to verify this system's capability to its full extent is given by curvilinear melting, which may be accompanied by serious complications when tested on a glacier. During the first EnEx project, there were only a few melting tests of this type in which the experience was gained that the IceMole could get stuck when being pulled backwards out of a curved melting channel.
Therefore, we tested the maneuverability of the probe  only for a relatively short-distance melting process of approximately three meters. \par
\tikzstyle{block} = [draw, fill=white, rectangle, 
    minimum height=3em, minimum width=4em]
\tikzstyle{sum} = [draw, fill=white, circle, node distance=1cm]
\tikzstyle{input} = [coordinate]
\tikzstyle{output} = [coordinate]
\tikzstyle{pinstyle} = [pin edge={to-,thin,black}]
\begin{figure}[h]
\centering
\vspace{.05cm}
\begin{tikzpicture}[auto, node distance=3cm,>=latex']
\begin{scope}[node distance=1.0cm and 10mm,every node/.style={scale=0.75}]
    \node [block, name=ocp] {TPP};
    \node [block, right= 1.5cm of ocp] (mpc) {MPC};
    \node [block, right= 1.5cm of mpc] (bct) {BCT};
    \node [block, below= of mpc] (fusion) {Sensor Fusion};
    \node [block, right= .5cm of fusion] (plant) {Plant};
    \node [block, right= 1.5cm of plant] (application) {SSG};
    \node[draw,dashed,scale=1.3,inner sep=3mm,label=below:IceMole,fit=(plant) (plant) (application) (plant)] {};
\end{scope}
    
    \draw [->] (ocp) -- node[name=u0] {\footnotesize $ \boldsymbol{x_{\text{init}}},\boldsymbol{x}_f$} (mpc);
    \draw [->] (mpc) -- node[name=ubar] {\footnotesize $\{\boldsymbol{u}^{[j]}\}^{[k]}$} (bct);
    \draw [->] (bct) -| node[name=vbar, near start] {\footnotesize $\{\boldsymbol{v}^{[j]}\}^{[k]}$} (application);
    \draw [->] (application) -- node[name=vbar] {\footnotesize $\{\boldsymbol{w}^{[j]}\}^{[k]}$} (plant);
    \draw [dashed,->] (plant) -- node[name=plantapp] {} (fusion);
    \draw [->] (fusion) -- node[name=fusiondata] {\footnotesize $\boldsymbol{x}^{[k]}$} (mpc);
    
\end{tikzpicture}
\caption{Process flow: In TPP a trajectory towards a certain destination $\boldsymbol{x}_f$ is planned. The state and control variables are used as an initial guess for MPC that frequently adjusts the trajectory with respect to the latest sensor fusion updates. A computed set of continuous control configurations $\{\boldsymbol{u}^{[j]}\}^{[k]}$ is forwarded to BCT, which recomputes a set of binary control configurations $\{\boldsymbol{v}^{[j]}\}^{[k]}$. Those are consecutively transmitted towards the IceMole interface where a software safeguard mechanism SSG has to be passed. A possibly different set of controls $\{\boldsymbol{w}^{[j]}\}^{[k]} \le \{\boldsymbol{v}^{[j]}\}^{[k]}$ is applied.}
\label{fig:methodology}
\end{figure}
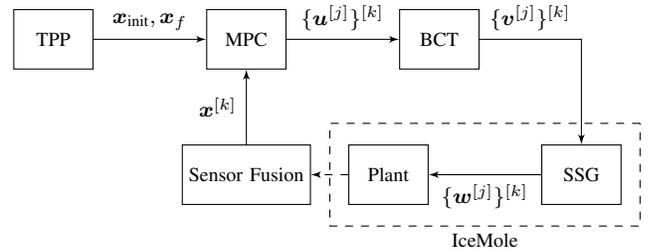
Before evaluation, we have to state to which terms the results could be obtained. The previously presented algorithms were performed in a sequence, as illustrated in \mbox{Fig. \ref{fig:methodology}}. Two things have not been discussed before. At first, there is a software safeguard mechanism (SSG) that might prevent the IceMole from switching certain heater units on, e.g. due to overheating. Secondly, we did not have any ground truth information with respect to the position data. Data from different subsystems, i.e. the IMU (inertial measurement unit), the EnEx-RANGE APU (Autonomous Pinger Units) network \cite{eliseev2014acoustic} and the screw rotation counter, was gathered and continuously fused to a pose by the sensor fusion subsystem \cite{clemens2017multi}. 
\par
To analyze the position data, we make the following definition. From now on we consider trajectories from $t_0$ to $t_k$ as a sequence of states
\begin{equation}
\label{eq:trajectorydefinition}
\mathrm{X}^{[k]} = [\boldsymbol{x}^{[0]},\dots,\boldsymbol{x}^{[k]}].
\end{equation}
 In general, we can choose a representation defined by a parameter setting $\boldsymbol{p}$ and a sequence of control variable sets $\boldsymbol{u}^{[0]},\dots,\boldsymbol{u}^{[k]}$.
We can calculate a trajectory ${{\mathrm{X}}_{\boldsymbol{u},\boldsymbol{p}}^{[k]} = [\boldsymbol{x}^{[0]},\dots,\tilde{\boldsymbol{x}}^{[k]}]}$ by forward integration by
\begin{equation}
\label{eq:forwardintegration}
\begin{array}{rl}
\tilde{\boldsymbol{x}}^{[k+1]} &= \tilde{\boldsymbol{x}}^{[k]}+\int\limits_{t_k}^{t_{k+1}}  f(\tilde{\boldsymbol{x}}^{[k]},\boldsymbol{u}^{[k]},\boldsymbol{p}) dt, \\
\tilde{\boldsymbol{x}}^{[0]}&=\boldsymbol{x}^{[0]}.
\end{array}
\end{equation}
For comparison, trajectories are generated by applying the integration scheme via the implicit trapezoidal rule to the different control set sequences, which have been computed and applied during the test run (cf. Fig. \ref{fig:methodology}). The parameter selection $\bar{\boldsymbol{p}}$ was fixed 
to values ${v_{max} = 1.0}$, ${\eta = 0.1}$, ${R_{\min} = 6.0}$ and ${\alpha = 0.3}$.  

\begin{figure}[h]
{
\centering
\Large
\includestandalone[width = 1\columnwidth]{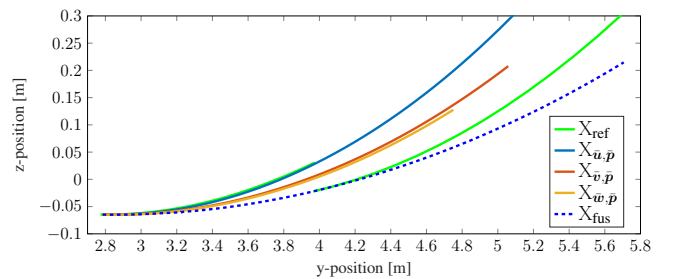}}
\vspace{-.25cm}
\caption[Comparison of trajectories]{Comparison of trajectories related to the different types of control sets. The two-part reference trajectory $\mathrm{X}_{\text{ref}}^{}$ diverges from the sensor fusion position data $\mathrm{X}_{\text{fus}}^{}$\footnotemark. Continuous ($\mathrm{X}_{\bar{\boldsymbol{u}},\bar{\boldsymbol{p}}}$), binary ($\mathrm{X}_{\bar{\boldsymbol{v}},\bar{\boldsymbol{p}}}$) and applied control sets ($\mathrm{X}_{\bar{\boldsymbol{w}},\bar{\boldsymbol{p}}}$) lie above as well.
The difference between $\mathrm{X}_{\bar{\boldsymbol{u}},\bar{\boldsymbol{p}}}$ and $\mathrm{X}_{\bar{\boldsymbol{w}},\bar{\boldsymbol{p}}}$ is significant, since in an ideal case, controls that come out of a model predictive control algorithm should be directly applicable.
}
\label{fig:trajectoriesAll}
\end{figure}
As can be deduced from Fig. \ref{fig:trajectoriesAll}, the IceMole has not been able to follow its planned trajectory $\mathrm{X}_{\text{ref}}^{}$\footnotetext{During the test, there was no information available regarding the yaw angle. The reference trajectory is two-part because execution had to be stopped due to reasons that are irrelevant for our methods and their analysis.} very well. All trajectories that have been simulated predict the IceMole to move in a much narrower curve than in reality. 
We explain these distinctions by two major reasons:
\begin{enumerate}
	\item approximation error introduced by the binary control transformation
	\item inaccuracy of the system dynamics model
\end{enumerate}
Below, we will locate these causes of error in detail and, consequently, improve our algorithms for future field tests.

\section{Post-Event Analysis}
\label{sec:analysis}
\subsection{Binary Control Transformation}
With the BCT algorithm, a real-time capable strategy to determine and alternate 32 binary controls has been postulated. Since the algorithm does not necessarily terminate in time, optimality cannot be expected. More importantly, the system dynamics are neglected. 
\begin{figure}[ht]
{
\centering
\Large
\includestandalone[width = 1\columnwidth]{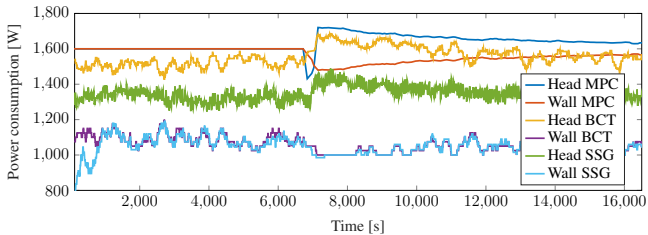}}
\caption{Simulated power consumption for head and wall heater units during the test (averaged over 1440 seconds). In case of the head heaters, the BCT algorithm manages to fit the continuous controls relatively well. There is a loss during application, though. The wall heater elements approximation is involved with a distinct deviation but the controls are nearly perfectly operated.}
\label{fig:powerAll}
\end{figure}
\par
\begin{figure}[h]
{
\centering
\Large
\vspace{.15cm}
\includestandalone[width = 1\columnwidth]{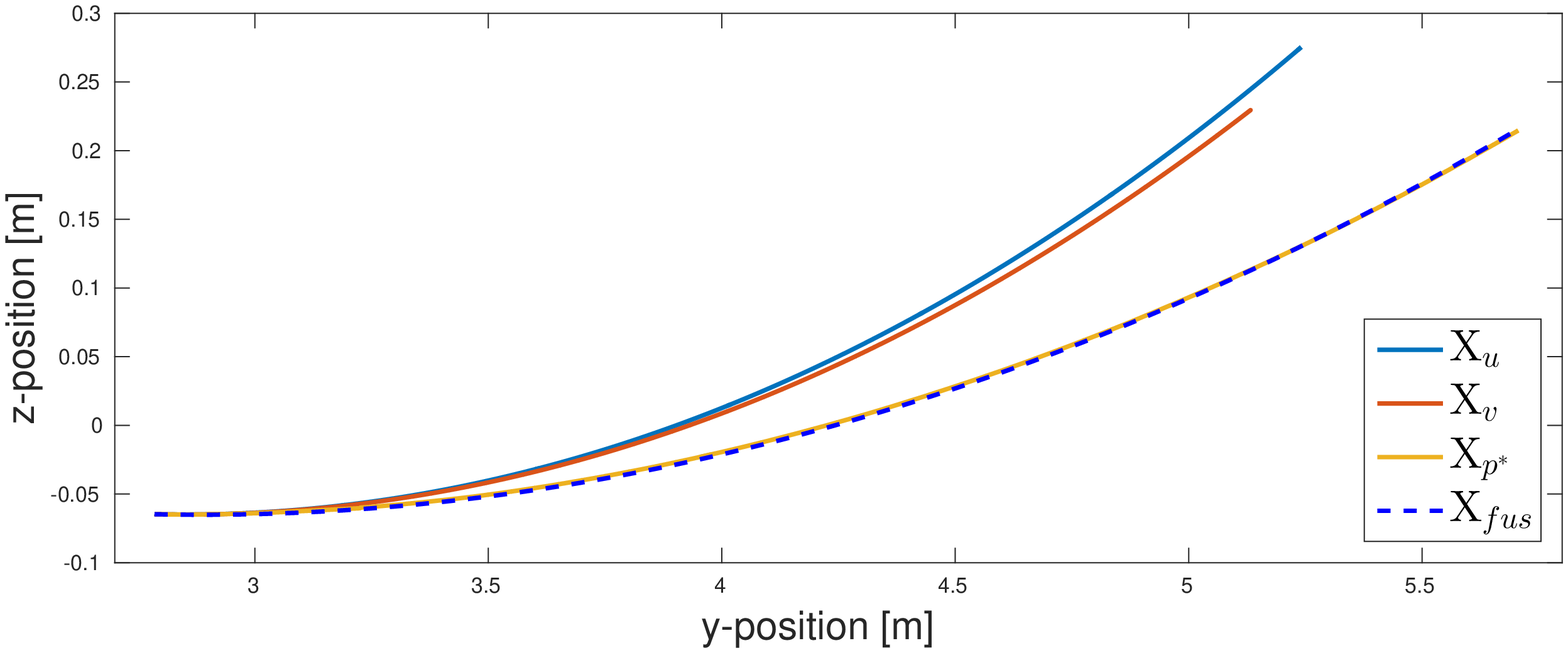}}
\caption{Trajectories related to modification of BCT algorithm and PI result. The deviation between the MPC result $\mathrm{X}_{\boldsymbol{u^{\ast}},\bar{\boldsymbol{p}}}$ and the BCT$^{\ast}$ result $\mathrm{X}_{\boldsymbol{v^{\ast}},\bar{\boldsymbol{p}}}$ is heavily reduced. With an optimal parameter set the respective trajectory $\mathrm{X}_{\bar{\boldsymbol{w}},\boldsymbol{p}^{\ast}}$ fits the fusion data $\mathrm{X}_{\text{fus}}$ perfectly.   }
	\label{fig:trajIdentBCT}
\end{figure}

 As can be deduced from Fig. \ref{fig:powerAll}, there is a systematic deviation regarding power distribution. 
During the test, we used weights ${\mu_{[i]}=1}$ for all controls within the objective function \eqref{eq:approximation}. A physically more reasonable choice would be to set ${\mu_{[i]}=P_{[i]}}$. This would lead to an optimization of the heating power distribution. Besides, since we are solving a fitting problem with additional constraints, the method tends to converge to a solution from below. To enforce a higher power output to the detriment of the uniformity of the distribution, we add a second term with a comparatively small weight ${\zeta > 0}$ to our cost function \eqref{eq:approximation}:
\begin{equation} 
\label{eq:unusedpower}
-\zeta \sum\limits_{i \in \mathcal{I}_T}^{}P_{[i]}\int_{0}^{T}v_{[i]}(t)\,dt
\end{equation}
Given these modifications, we can define the modified binary control transformation by
\begin{equation}
\label{eq:bct*}
\tag{BCT$^{\ast}$}
\hspace*{-.30cm}
\begin{array}{lll}
\min\limits_{v} &\sum\limits_{i \in \mathcal{I}_T}^{}P_{[i]}^2\sum\limits_{j=1}^{M} a^{[j]} (\sum\limits_{m=j-l}^{j} b^{[m]} (u_{[i]}^{[m]}-v_{[i]}^{[m]}) )^2  \\[-5pt]
& - \zeta\sum\limits_{i \in \mathcal{I}_T}^{}P_{[i]} \sum\limits_{j=1}^{M} a^{[j]}v_{[i]}^{[j]} \\[10pt]
s.t. & \text{eq. }\eqref{eq:powerall} \text{ hold }\, \forall j = 1,\dots,M.
\end{array}
\end{equation}
Additionally, there was an inconsistency regarding the power restrictions \eqref{eq:powerallgroup} and \eqref{eq:powerwallgroup} between MPC and BCT, which has been eliminated.
We performed a reoptimization with the modified algorithms and got a much better result, as can be seen in Fig. \ref{tab:power} with respect to the power distribution, and in Fig. \ref{fig:trajIdentBCT} with respect to the trajectory.
\begin{figure}[h]
\centering
\begin{tabular}{cp{1.2cm}p{1.2cm}p{1.2cm}}
\toprule
\multirow{ 2}{*}{} & \multicolumn{3}{c @{}}{Power Consumption (rel. to MPC)} \\ \cmidrule(l){2-4}  
 & BCT & SSG & BCT*\\   \midrule
Head ($\mathcal{I}_H$) & 0.96 & 0.82 & 0.95\\    \hline
Wall ($\mathcal{I}_W$)& 0.68 & 0.67 & 0.98\\    \hline
Back ($\mathcal{I}_B$)& 0.99 & 0.80 & 0.95\\ 
\bottomrule
\end{tabular}
\caption{Power consumption of different binary control sets relative to MPC controls $\boldsymbol{\bar u}$. During the test, the BCT algorithm generated a set of controls $\boldsymbol{\bar v}$ in which wall heater units were heavily underrepresented. The software safeguard mechanism SSG prohibited a direct application $\boldsymbol{\bar w}$ of head and back heater units though. With the modified algorithm BCT$^{\ast}$ the coverage is more than 95 percent for all groups.}
\label{tab:power}
\end{figure}
\subsection{Parameter Identification}
To a certain extent we have to accept that our model is not precise enough. We neither take any heating flow concepts into account nor kinetic equations. 
Additionally, we do not consider locally changing environmental conditions, e.g. small crevasses. 
Overall, we cannot expect physical correctness for a highly complex melting process if a model is that simple. Nevertheless, it should be of higher precision for this short-distance test case. \par
To overcome the inaccuracies, we will not extend our model in the first place, as its simplicity allows complex trajectory planning. Instead, we will adapt our system parameters. To recalculate a better choice, we state another optimization problem.
With the forward integration scheme proposed in \eqref{eq:forwardintegration},
a trajectory ${\mathrm{X}}_{\bar{\boldsymbol{w}},\boldsymbol{p}}$ as a result of a known sequence of applied binary controls can be calculated with respect to parameters $\boldsymbol{p}$.
To fit the fusion position data, we define a cost function by the least squares error with respect to the y- and z-position:
\begin{equation}
F({\mathrm{X}}_{\bar{\boldsymbol{w}},\boldsymbol{p}},\mathrm{X}_{\text{fus}})
:= \sum\limits_{l=1}^{N} ({y_{\bar{\boldsymbol{w}},\boldsymbol{p}}}^{[l]}-y_{\text{fus}}^{[l]})^2+({z_{\bar{\boldsymbol{w}},\boldsymbol{p}}}^{[l]}-z_{\text{fus}}^{[l]})^2
\end{equation}
Finally, the parameter identification can be stated as an unconstrained optimization problem
\begin{equation}
\label{eq:pi}
\tag{PI}
\min\limits_{\boldsymbol{p}}  F({\mathrm{X}}_{\bar{\boldsymbol{w}},\boldsymbol{p}},\mathrm{X}_{\text{fus}}).
\end{equation}
We solve PI and obtain an optimal parameter set ${\boldsymbol{p}^{\ast} = \arg\min_{\boldsymbol{p}} F({\mathrm{X}}_{\bar{\boldsymbol{w}},\boldsymbol{p}},\mathrm{X}_{\text{fus}})}$ with ${v_{max} = 1.48}$, ${\eta = 0.10}$, ${R_{\min} = 8.66}$ and ${\alpha = 0.22}$.
Except for $v_{max}$, this parameter set is reasonable with respect to its physical background. The comparatively high value of $v_{max}$ (around +20$\%$) can partially be explained by the effective loss of induced heating power along the process chain.
\subsection{Combination}
By modifying the BCT algorithm and using parameter identification, we resolved the different sources of error separately.
However, on the basis of the given data set, we cannot bring these parts together without deliberating. 
The challenge would be to find a parameter set $\boldsymbol{p^{\ast}}$ that is feasible at planning time. Because of the step-wise decrease of induced heating power along the process sequence, the previously computed parameter set $\boldsymbol{p^{\ast}}$ is not directly applicable. To overcome this, there are two possible approaches. One approach would be to integrate BCT$^{\ast}$ and SSG into the MPC algorithm. As to BCT$^{\ast}$, it is too computationally expensive and contradicts the separation of a MIOCP in the first place. SSG is not an accessible algorithm. Therefore, both would have to be approximated with respect to their effects. Alternatively, one could try to ``invert'' the applied controls $\boldsymbol{\bar w}$ at first, in order to perform PI for a continuous control set $\boldsymbol{u}$. Again, this requires an approximation of the inverse operator. Both strategies could be carried out by introducing parameter-dependent functions that model the assumed decrease of induced heating power along the sequence of algorithms. 

\section{Conclusions \& Future Work}
\label{sec:conclusion}
We presented an optimization-based framework for the maneuvering of an autonomous melting probe. It includes a strategy for reformulating and separating a MIOCP into two consecutive algorithms, so that an altogether strong solution can be found under real-time conditions. The framework's fundamental applicability was proven during an intermediate glacier test. An evaluation of the data indicated systematic errors which could be attributed to different sources. Those defects were fixed separately by using parameter identification and adapting an objective function. A combination of these strategies yields a new problem which will become a subject to further investigations.
\par As a result, we have to adapt our framework so that the cycle between trajectory (re-)planning and parameter identification can be closed. We will improve our models with the objective to integrate an on-line adaptive model for the next glacier test.

\addtolength{\textheight}{-12cm}   




\section*{ACKNOWLEDGMENT}

This work was supported by the German Aerospace Center
(DLR) with financial means of the German Federal Ministry
for Economic Affairs and Energy (BMWi), project ``EnEx-CAUSE''
(grant No. 50 NA 1505).

\bibliographystyle{IEEEtran}
\bibliography{IEEEabrv,references}

\begin{thebibliography}{10}
\providecommand{\url}[1]{#1}
\csname url@samestyle\endcsname
\providecommand{\newblock}{\relax}
\providecommand{\bibinfo}[2]{#2}
\providecommand{\BIBentrySTDinterwordspacing}{\spaceskip=0pt\relax}
\providecommand{\BIBentryALTinterwordstretchfactor}{4}
\providecommand{\BIBentryALTinterwordspacing}{\spaceskip=\fontdimen2\font plus
\BIBentryALTinterwordstretchfactor\fontdimen3\font minus
  \fontdimen4\font\relax}
\providecommand{\BIBforeignlanguage}[2]{{%
\expandafter\ifx\csname l@#1\endcsname\relax
\typeout{** WARNING: IEEEtran.bst: No hyphenation pattern has been}%
\typeout{** loaded for the language `#1'. Using the pattern for}%
\typeout{** the default language instead.}%
\else
\language=\csname l@#1\endcsname
\fi
#2}}
\providecommand{\BIBdecl}{\relax}
\BIBdecl

\bibitem{porco2006cassini}
C.~Porco, P.~Helfenstein, P.~Thomas, A.~Ingersoll, J.~Wisdom, R.~West,
  G.~Neukum, T.~Denk, R.~Wagner, T.~Roatsch \emph{et~al.}, ``{Cassini Observes
  the Active South Pole of Enceladus},'' \emph{Science}, vol. 311, no. 5766,
  pp. 1393--1401, 2006.

\bibitem{dachwald2014icemole}
B.~Dachwald, J.~Mikucki, S.~Tulaczyk, I.~Digel, C.~Espe, M.~Feldmann,
  G.~Francke, J.~Kowalski, and C.~Xu, ``{IceMole: a maneuverable probe for
  clean in situ analysis and sampling of subsurface ice and subglacial aquatic
  ecosystems},'' \emph{Annals of Glaciology}, vol.~55, no.~65, pp. 14--22,
  2014.

\bibitem{kowalski2016navigation}
J.~Kowalski, P.~Linder, S.~Zierke, B.~von Wulfen, J.~Clemens,
  K.~Konstantinidis, G.~Ameres, R.~Hoffmann, J.~Mikucki, S.~Tulaczyk
  \emph{et~al.}, ``{Navigation Technology for Exploration of Glacier Ice With
  Maneuverable Melting Probes},'' \emph{Cold Regions Science and Technology},
  vol. 123, pp. 53--70, 2016.

\bibitem{sussmann1997300}
H.~J. Sussmann and J.~C. Willems, ``{300 Years of Optimal Control: From the
  Brachystochrone to the Maximum Principle},'' \emph{IEEE Control Systems
  Magazine}, vol.~17, no.~3, pp. 32--44, 1997.

\bibitem{stryk1992direct}
O.~Von~Stryk and R.~Bulirsch, ``Direct and indirect methods for trajectory
  optimization,'' \emph{{Annals of Operations Research}}, vol.~37, no.~1, pp.
  357--373, 1992.

\bibitem{kirches2011fast}
C.~Kirches, \emph{Fast Numerical Methods for Mixed-Integer Nonlinear
  Model-Predictive Control}.\hskip 1em plus 0.5em minus 0.4em\relax Springer,
  2011.

\bibitem{gerdts2006variable}
M.~Gerdts, ``A variable time transformation method for mixed-integer optimal
  control problems,'' \emph{Optimal Control Applications and Methods}, vol.~27,
  no.~3, pp. 169--182, 2006.

\bibitem{sager2005numerical}
S.~Sager, \emph{Numerical methods for mixed-integer optimal control
  problems}.\hskip 1em plus 0.5em minus 0.4em\relax Der andere Verlag
  T{\"o}nning, L{\"u}beck, Marburg, 2005.

\bibitem{aamot1967heat}
H.~Aamot, \emph{Heat Transfer and Performance Analysis of a Thermal Probe for
  Glaciers}, ser. Technical report.\hskip 1em plus 0.5em minus 0.4em\relax U.S.
  Army Materiel Command, Cold Regions Research \& Engineering Laboratory, 1967.

\bibitem{ulamec2007access}
S.~Ulamec, J.~Biele, O.~Funke, and M.~Engelhardt, ``Access to glacial and
  subglacial environments in the solar system by melting probe technology,''
  \emph{Reviews in Environmental Science and Bio/Technology}, vol.~6, no. 1-3,
  pp. 71--94, 2007.

\bibitem{kumano2005direct}
H.~Kumano, A.~Saito, S.~Okawa, and Y.~Yamada, ``Direct contact melting with
  asymmetric load,'' \emph{International Journal of Heat and Mass Transfer},
  vol.~48, no.~15, pp. 3221--3230, 2005.

\bibitem{schuller2016curvilinear}
K.~Sch{\"u}ller, J.~Kowalski, and P.~R{\aa}back, ``Curvilinear melting--a
  preliminary experimental and numerical study,'' \emph{International Journal
  of Heat and Mass Transfer}, vol.~92, pp. 884--892, 2016.

\bibitem{stengel2015flight}
R.~F. Stengel, \emph{Flight dynamics}.\hskip 1em plus 0.5em minus 0.4em\relax
  Princeton University Press, 2015.

\bibitem{kavraki1996probabilistic}
L.~E. Kavraki, P.~Svestka, J.-C. Latombe, and M.~H. Overmars, ``Probabilistic
  roadmaps for path planning in high-dimensional configuration spaces,''
  \emph{IEEE transactions on Robotics and Automation}, vol.~12, no.~4, pp.
  566--580, 1996.

\bibitem{lavalle98rapidly}
S.~M. Lavalle, ``Rapidly-exploring random trees: A new tool for path
  planning,'' Tech. Rep., 1998.

\bibitem{bixby2004mixed}
R.~E. Bixby, M.~Fenelon, Z.~Gu, E.~Rothberg, and R.~Wunderling, ``Mixed integer
  programming: a progress report,'' \emph{The Sharpest Cut}, pp. 309--325,
  2004.

\bibitem{buskens2012esa}
C.~B{\"u}skens and D.~Wassel, ``{The ESA NLP Solver Worhp},'' in \emph{Modeling
  and Optimization in Space Engineering}.\hskip 1em plus 0.5em minus
  0.4em\relax Springer, 2012, pp. 85--110.

\bibitem{bloesch2016primer}
M.~Bloesch, H.~Sommer, T.~Laidlow, M.~Burri, G.~Nuetzi, P.~Fankhauser,
  D.~Bellicoso, C.~Gehring, S.~Leutenegger, M.~Hutter \emph{et~al.}, ``A primer
  on the differential calculus of 3d orientations,'' \emph{arXiv preprint
  arXiv:1606.05285}, 2016.

\bibitem{grune2011nonlinear}
L.~Gr{\"u}ne and J.~Pannek, ``Nonlinear model predictive control,'' in
  \emph{Nonlinear Model Predictive Control}.\hskip 1em plus 0.5em minus
  0.4em\relax Springer, 2011, pp. 43--66.

\bibitem{scip2017}
S.~J. Maher, T.~Fischer, T.~Gally, G.~Gamrath, A.~Gleixner, R.~L. Gottwald,
  G.~Hendel, T.~Koch, M.~E. L{\"u}bbecke, M.~Miltenberger, B.~M{\"u}ller, M.~E.
  Pfetsch, C.~Puchert, D.~Rehfeldt, S.~Schenker, R.~Schwarz, F.~Serrano,
  Y.~Shinano, D.~Weninger, J.~T. Witt, and J.~Witzig,
  ``\BIBforeignlanguage{eng}{{The SCIP Optimization Suite 4.0}},'' ZIB,
  Takustr.7, 14195 Berlin, Tech. Rep. 17-12, 2017.

\bibitem{eliseev2014acoustic}
D.~Eliseev, D.~Heinen, K.~Helbing, R.~Hoffmann, U.~Naumann, F.~Scholz,
  C.~Wiebusch, and S.~Zierke, ``Acoustic in-ice positioning in the enceladus
  explorer project,'' \emph{Annals of Glaciology}, vol.~55, no.~68, pp.
  253--259, 2014.

\bibitem{clemens2017multi}
J.~Clemens, ``Multi-robot in-ice localization using graph optimization,'' in
  \emph{Autonomous Robot Systems and Competitions (ICARSC), 2017 IEEE
  International Conference on}.\hskip 1em plus 0.5em minus 0.4em\relax IEEE,
  2017, pp. 36--42.

\end{thebibliography}

\end{document}